\documentclass[11pt]{article}
\usepackage{amsmath, amsthm}
\usepackage{url}

\newtheoremstyle{theorem}
  {10pt}		  
  {10pt}  
  {\sl}  
  {\parindent}     
  {\bf}  
  {. }    
  { }    
  {}     
\theoremstyle{theorem}
\newtheorem{theorem}{Theorem}
\newtheorem{corollary}[theorem]{Corollary}
\newtheorem{lemma}[theorem]{Lemma}
\newtheorem{conjecture}[theorem]{Conjecture}

\newtheoremstyle{defi}
  {10pt}		  
  {10pt}  
  {\rm}  
  {\parindent}     
  {\bf}  
  {. }    
  { }    
  {}     
\theoremstyle{defi}

\newtheorem{remark}[theorem]{Remark}

\begin{document}

\title{New Results for the Descartes-Frenicle-Sorli \\
Conjecture on Odd Perfect Numbers}

\author{Jose Arnaldo B. Dris\\
Far Eastern University\\
Manila, Philippines\\
josearnaldobdris@gmail.com\\[2pt]}

\maketitle

\begin{abstract}
If $N={q^k}{n^2}$ is an odd perfect number given in Eulerian form, then the Descartes-Frenicle-Sorli conjecture predicts that $k=1$. Brown \cite{Brown} has recently announced a proof for the inequality \\$q < n$, and a partial proof that $q^k < n$ holds under many cases. In this article, we give a strategy for strengthening Brown's result to $q^2 < n$.

{\bf AMS Subject Classification: Primary 11A05; Secondary 11J25, 11J99}

{\bf Key Words and Phrases: odd perfect number, Sorli's conjecture, Euler prime}
\end{abstract}

\section{Introduction}\label{Section1}
If $N$ is a positive integer, then we write $\sigma(N)$ for the sum of the divisors of $N$.  A number $N$ is \emph{perfect} if $\sigma(N)=2N$.  It is currently unknown whether there are infinitely many even perfect numbers, or whether any odd perfect numbers (OPNs) exist.  Ochem and Rao recently proved \cite{OchemRao} that, if $N$ is an odd perfect number, then $N > {10}^{1500}$ and that the largest component (i.e., divisor $p^a$ with $p$ prime) of $N$ is bigger than ${10}^{62}$.  This improves on previous results by Brent, Cohen and te Riele \cite{BrentCohenteRiele} in 1991 ($N > {10}^{300}$) and Cohen \cite{Cohen} in 1987 (largest component $p^a > {10}^{20}$). 

An odd perfect number $N = {q^k}{n^2}$ is said to be given in Eulerian form if $q$ is prime with $q \equiv k \equiv 1 \pmod 4$ and $\gcd(q, n) = 1$.  (The number $q$ is called the \emph{Euler prime}, while the component $q^k$ is referred to as the \emph{Euler factor}.  Note that, since $q$ is prime and $q \equiv 1 \pmod 4$, then $q \geq 5$.)

We denote the abundancy index $I$ of the positive integer $x$ as $$I(x) = \frac{\sigma(x)}{x}.$$

In his Ph.~D.~thesis, Sorli \cite{Sorli} conjectured that $k=1$, after testing large numbers with $8$ distinct prime factors for perfection. (More recently, Beasley \cite{Beasley} points out that Descartes was the first to conjecture $k=1$ ``in a letter to Mersenne in 1638, with Frenicle's subsequent observation occurring in 1657".)

In the M.~Sc.~thesis \cite{Dris3}, it was conjectured that the components $q^k$ and $n$ are related by the inequality $q^k < n$.  This conjecture was made on the basis of the result $I(q^k)<I(n)$.  Recently, Brown \cite{Brown} announced a proof for the inequality $q < n$, and a partial proof that $q^k < n$ holds under many cases.

\section{Conditions Sufficient for Sorli's Conjecture}\label{Section2}
Some sufficient conditions for Sorli's conjecture were given in \cite{Dris}.  We reproduce these conditions here.

\begin{lemma}\label{Lemma1}
Let $N = {q^k}{n^2}$ be an odd perfect number given in Eulerian form.  If $n<q$, then $k=1$.
\end{lemma}

\begin{remark}\label{Remark1}
The proof of Lemma \ref{Lemma1} follows from the inequality $q^k < n^2$ and the congruence $k \equiv 1 \pmod 4$ (see \cite{Dris}).  (Note the related inequality $$I(q^k)<I(n^2)$$
for the abundancy indices of the components $q^k$ and $n^2$.)
\end{remark}

\begin{lemma}\label{Lemma2}
Let $N = {q^k}{n^2}$ be an odd perfect number given in Eulerian form.  If 
$$\sigma(n) \le \sigma(q),$$ 
then $k=1$.
\end{lemma}

\begin{lemma}\label{Lemma3}
Let $N = {q^k}{n^2}$ be an odd perfect number given in Eulerian form.  If $$\frac{\sigma(n)}{q} < \frac{\sigma(q)}{n},$$ 
then $k=1$.
\end{lemma}

\begin{remark}\label{Remark2}
Notice that, if $$\frac{\sigma(n)}{q} < \frac{\sigma(q)}{n},$$
then it follows that $$\frac{\sigma(n)}{q^k} = \frac{\sigma(n)}{q} < \frac{\sigma(q)}{n} = \frac{\sigma(q^k)}{n}.$$  
Consequently, by the contrapositive, if $$\frac{\sigma(q^k)}{n} < \frac{\sigma(n)}{q^k},$$
then $$\frac{\sigma(q)}{n} \leq \frac{\sigma(q^k)}{n} < \frac{\sigma(n)}{q^k} \leq \frac{\sigma(n)}{q}.$$ 
\end{remark}

\begin{remark}\label{Remark3}
Let $N = {q^k}{n^2}$ be an odd perfect number given in Eulerian form. Suppose that $$\frac{\sigma(q)}{n} = \frac{\sigma(n)}{q}.$$  
Then we know that:
$${q}{\sigma(q)} = {n}{\sigma(n)}.$$

Since $\gcd(q, n) = 1$, then $q \mid \sigma(n)$ and $n \mid \sigma(q)$.  Therefore, it follows that $\displaystyle\frac{\sigma(q)}{n}$ and $\displaystyle\frac{\sigma(n)}{q}$ are equal positive integers.

This is a contradiction, as:
$$1 < I(q) = \frac{\sigma(q)}{q} = 1 + \frac{1}{q} \leq \frac{6}{5} < \sqrt{\frac{5}{3}} < I(n) < I(q)I(n) = I(qn) < 2$$
which implies that:
$$1 < \sqrt{\frac{5}{3}} < I(n) < I(q)I(n) = I(qn)
= \left[\frac{\sigma(q)}{q}\right]\left[\frac{\sigma(n)}{n}\right]
= \left[\frac{\sigma(q)}{n}\right]\left[\frac{\sigma(n)}{q}\right] < 2$$

Consequently,
$$\frac{\sigma(q)}{n} \neq \frac{\sigma(n)}{q}.$$
Similarly, we can prove that
$$\frac{\sigma(q^k)}{n} \neq \frac{\sigma(n)}{q^k}.$$
\end{remark}

\begin{lemma}\label{Lemma4}
Let $N = {q^k}{n^2}$ be an odd perfect number given in Eulerian form.  Then $n<q$ if and only if $N<{q^3}$.
\end{lemma}

\begin{proof}
Suppose that $N = {q^k}{n^2}$ is an odd perfect number given in Eulerian form. If $n<q$, then assuming to the contrary that ${q^3}<N$, we get that
$${q^3}<N=q{n^2}<q\cdot{q^2}=q^3$$
since $n<q$ implies $k=1$, by Lemma \ref{Lemma1}.
For the other direction, if $N<{q^3}$, then ${q^k}{n^2}<{q^3}$, so that we have
$${n^2}<{q^{3-k}}\leq{q^2}$$
since $k \equiv 1 \pmod 4$ implies that $k \geq 1$. Consequently, $n < q$, and we are done.
\end{proof}

\begin{corollary}\label{CorollaryToLemma4}
Let $N = {q^k}{n^2}$ be an odd perfect number given in Eulerian form.  Then $n<q^{5/2}$ if and only if $N<{q^6}$.
\end{corollary}

\begin{proof}
First we show that $n<q^{5/2}$ implies $k=1$. To this end, assuming $n<q^{5/2}$, since $q^k < n^2$ (see \cite{Dris}), we then have that:
$$q \leq q^k < n^2 < q^5.$$
The last chain of inequalities implies that
$$1 \leq k < 5.$$
This inequality, together with the condition $k \equiv 1 \pmod 4$, implies that $k=1$.

We now prove the claim in Corollary \ref{CorollaryToLemma4}. If $n<q^{5/2}$, then assuming to the contrary that ${q^6}<N$, we get that
$${q^6}<N=q{n^2}<q\cdot{q^5}={q^6}.$$
This is a contradiction.
For the other direction, if $N<{q^6}$, then ${q^k}{n^2}<{q^6}$, so that we have
$${n^2}<{q^{6-k}}\leq{q^5}$$
since $k \equiv 1 \pmod 4$ implies that $k \geq 1$. Consequently, $n < q^{5/2}$, and we are done.
\end{proof}

\begin{remark}\label{Remark4}
A recent result by Acquaah and Konyagin \cite{AcquaahKonyagin} \emph{almost} disproves $n<q$. They obtained the estimate $y < (3N)^{1/3}$ for all the prime factors $y$ of an odd perfect number $N$.  In particular, if $N={q^k}{n^2}$ is an odd perfect number given in Eulerian form, then letting $y=q$ and assuming $k=1$ gives:
$$q < (3N)^{1/3} = (3q{n^2})^{1/3} \Longrightarrow q^3 < 3q{n^2} \Longrightarrow q < n\sqrt{3}.$$

Since the contrapositive of the implication $n<q \Longrightarrow k=1$ is $k > 1 \Longrightarrow q < n$, it follows that the inequality
$$q<n\sqrt{3}$$
holds unconditionally, regardless of the status of Sorli's conjecture.

More recently, Brown \cite{Brown} claims a proof for the inequality $q < n$, and a partial proof that $q^k < n$ holds under many cases.
\end{remark}

We now give a condition that is weaker than $n < q$, which also implies $k=1$.

\begin{lemma}\label{WeakConditionImpliesSorli}
Let $N = {q^k}{n^2}$ be an odd perfect number given in Eulerian form. Then
$$n < {\left(\frac{3}{2}{q^5}\right)}^{1/2}$$
implies $k=1$.
\end{lemma}

\begin{proof}
Suppose that $N = {q^k}{n^2}$ is an odd perfect number given in Eulerian form.  Let
$$n < {\left(\frac{3}{2}{q^5}\right)}^{1/2}$$
and assume to the contrary that $k \neq 1$.  Since $k \equiv 1 \pmod 4$, this means that $k \geq 5$.
Additionally, from \cite{Dris}, we have that
$$q^k < \sigma(q^k) \leq \frac{2}{3}{n^2}.$$
Consequently, we have the following chain of inequalities:
$$q^5 \leq q^k < \frac{2}{3}{\left({\left(\frac{3}{2}{q^5}\right)}^{1/2}\right)}^2 < q^5.$$
This is a contradiction.
\end{proof}

We also have the following corollary to Lemma \ref{WeakConditionImpliesSorli}, and this uses a result from \cite{BroughanDelbourgoZhou}.

\begin{corollary}\label{CorollarytoWeakConditionImpliesSorli}
Let $N = {q^k}{n^2}$ be an odd perfect number given in Eulerian form. Then
$$n < {\left(\frac{315}{2}{q^5}\right)}^{1/2}$$
implies $k=1$.
\end{corollary}

\begin{proof}
The proof is very similar to that of Lemma \ref{WeakConditionImpliesSorli}, except that it uses the improved bound
$$\sigma(q^k) \leq \frac{2}{315}{n^2}$$
(see \cite{BroughanDelbourgoZhou}) instead of
$$\sigma(q^k) \leq \frac{2}{3}{n^2}$$
(see \cite{Dris}).
\end{proof}

\begin{remark}\label{ConditionsEquivalentToLemma10AndCorollary11}
Similar to the proofs of Lemma \ref{Lemma4} and Corollary \ref{CorollaryToLemma4}, we can show that the following biconditionals are true:
$$n < {\left(\frac{3}{2}{q^5}\right)}^{1/2} \Longleftrightarrow N < \frac{3}{2}{q^6}$$
$$n < {\left(\frac{315}{2}{q^5}\right)}^{1/2} \Longleftrightarrow N < \frac{315}{2}{q^6}$$
\end{remark}

\begin{remark}\label{ChenChenImprovingOnBroughanDelbourgoZhou}
Chen and Chen \cite{ChenChen} has a relatively recent paper which further improves on Broughan et.~al.'s results (see \cite{BroughanDelbourgoZhou}).  They also pose a related open problem.
\end{remark}

\section{New Results Related to Sorli's Conjecture}\label{Section3}
First, we reproduce the following lemma from \cite{Dris}, as we will be using these results later.

\begin{lemma}\label{Lemma5}
Let $N = {q^k}{n^2}$ be an odd perfect number given in Eulerian form.  The following series of inequalities hold:
\begin{itemize}
\item{If $k = 1$, then $1 < I(q^k) = I(q) \leq \frac{6}{5} < \sqrt{\frac{5}{3}} < I(n) < 2$.}
\item{If $k \geq 1$, then $1 < I(q^k) < \frac{5}{4} < \sqrt{\frac{8}{5}} < I(n) < 2$.}
\end{itemize}
\end{lemma}

We have the following (slightly) stronger inequality from \cite{Dris}.

\begin{lemma}\label{Lemma6}
Let $N = {q^k}{n^2}$ be an odd perfect number given in Eulerian form.  Then ${\left(I(q^k)\right)}^2 < I(n^2)$.
\end{lemma}

\begin{proof}
The proof follows from the inequality $I(q^k) < \sqrt[3]{2}$ and the equation $2 = I(q^k)I(n^2)$.
\end{proof}

\begin{remark}\label{Remark5}
Another proof of Lemma \ref{Lemma6} is as follows:
$$I(q^k) < \frac{5}{4} \Longrightarrow {\left(I(q^k)\right)}^2 < \frac{25}{16} = 1.5625 < 1.6 = \frac{8}{5} < I(n^2).$$
In fact, if
$${\left(I(q^k)\right)}^y < {\left(\frac{5}{4}\right)}^y \leq \frac{8}{5} < I(n^2)$$
then
$$y \leq \frac{3\log{2} - \log{5}}{\log{5} - 2\log{2}}.$$
Thus, if we let 
$$z = \frac{3\log{2} - \log{5}}{\log{5} - 2\log{2}} \approx 2.1062837195,$$ 
then
$${\left(I(q^k)\right)}^z \leq \frac{8}{5} < I(n^2).$$
\end{remark}

Next, we derive a lower bound for $I(q^k)+I(n)$.

\begin{lemma}\label{Lemma7}
Let $N = {q^k}{n^2}$ be an odd perfect number given in Eulerian form.  The following inequality holds:
$$I(q^k) + I(n) \geq I(q) + I(n) > 1 + \sqrt{2}.$$  
\end{lemma}

\begin{proof}
Let $N = {q^k}{n^2}$ be an odd perfect number given in Eulerian form.  Then we have the following:
$$I(q^k) + I(n) \geq I(q) + I(n) \geq 1 + \frac{1}{q} + \sqrt{\frac{2(q-1)}{q}}.$$
But
$$f(q) = 1 + \frac{1}{q} + \sqrt{\frac{2(q-1)}{q}}$$
is a decreasing function of $q$.  Consequently,
$$f(q) > \lim_{q \rightarrow \infty}{\left(1 + \frac{1}{q} + \sqrt{\frac{2(q-1)}{q}}\right)} = 1 + \sqrt{2}.$$
\end{proof}

\begin{remark}\label{Remark6}
The following result was communicated to the author (via e-mail, by Pascal Ochem) in April of 2013: If $N={q^k}{n^2}$ is an odd perfect number given in Eulerian form, then
$$I(n)>{\left(\frac{8}{5}\right)}^{\frac{\ln(4/3)}{\ln(13/9)}} \approx 1.44440557.$$
(Note that ${\left(\frac{8}{5}\right)}^{\frac{\ln(4/3)}{\ln(13/9)}} > \sqrt{2}$.)
\end{remark}

Further to Remark \ref{Remark6} and Lemma \ref{Lemma6}, we have the following related result.

\begin{lemma}\label{CorollaryToLemma6}
Let $N = {q^k}{n^2}$ be an odd perfect number given in Eulerian form.  Then ${\left(I(q)\right)}^2 < I(n)$.
\end{lemma}

\begin{proof}
By Lemma \ref{Lemma5},
$$I(q) \leq \frac{6}{5} \Longrightarrow {\left(I(q)\right)}^2 \leq \frac{36}{25} = 1.44.$$
The conclusion follows from the result $I(n)>1.44440557$ in Remark \ref{Remark6}.

In fact, if
$$\left(I(q)\right)^u < \left(\frac{6}{5}\right)^u \leq {\left(\frac{8}{5}\right)}^{\frac{\ln(4/3)}{\ln(13/9)}}$$
then
$$u \leq -\frac{\left(2\log(2) - \log(3)\right)\left(3\log(2) - \log(5)\right)}{\left(\log(2) + \log(3) - \log(5)\right)\left(2\log(3) - \log(13)\right)}.$$
Thus, if we let
$$v = -\frac{\left(2\log(2) - \log(3)\right)\left(3\log(2) - \log(5)\right)}{\left(\log(2) + \log(3) - \log(5)\right)\left(2\log(3) - \log(13)\right)} \approx 2.0168$$
then
$$\left(I(q)\right)^v \leq {\left(\frac{8}{5}\right)}^{\frac{\ln(4/3)}{\ln(13/9)}} < I(n).$$
\end{proof}

\begin{remark}\label{Remark7}
As pointed out by Ochem to the author (via the same e-mail mentioned in Remark \ref{Remark6}), a case-by-case analysis yields a sharper lower bound for $I(q^k)+I(n)$:
\begin{itemize}
\item{If $q = 5$ then $I(q^k)+I(n) \geq I(q)+I(n) \geq (6/5) + {(8/5)}^{\ln(4/3)/\ln(13/9)} \approx 2.6444055$.}
\item{If $q \geq 13$ then $I(q^k)+I(n) \geq I(q)+I(n) \geq (14/13) + {(24/13)}^{\ln(4/3)/\ln(13/9)} \approx 2.6924318$.}
\end{itemize}
Therefore, we have the lower bound
$$I(q^k)+I(n) \geq I(q)+I(n) \geq \frac{6}{5} + {\left(\frac{8}{5}\right)}^{\frac{\ln(4/3)}{\ln(13/9)}} \approx 2.6444055.$$
\end{remark}

We now state and prove the following theorem, which provides conditions equivalent to the conjecture mentioned in the introduction.

\begin{theorem}\label{Theorem1}
If $N = {q^k}{n^2}$ is an odd perfect number given in Eulerian form, then the following biconditional is true:
$$q^k < n \Longleftrightarrow \sigma(q^k) < \sigma(n).$$
\end{theorem}

In preparation for the proof of Theorem \ref{Theorem1}, we derive the following results.

\begin{lemma}\label{Lemma8}
Let $N = {q^k}{n^2}$ be an odd perfect number given in Eulerian form.  If
$$I(q^k) + I(n) < \frac{\sigma(q^k)}{n} + \frac{\sigma(n)}{q^k},$$
then
$$q^k < n \Longleftrightarrow \sigma(q^k) < \sigma(n).$$
\end{lemma}

\begin{proof}
Let $N = {q^k}{n^2}$ be an odd perfect number given in Eulerian form.
Assume that $$I(q^k) + I(n) < \frac{\sigma(q^k)}{n} + \frac{\sigma(n)}{q^k}.$$
It follows that
$$I(q^k) + I(n) < \left(\frac{q^k}{n}\right)I(q^k) + \left(\frac{n}{q^k}\right)I(n).$$
Consequently,
$${{q^k}n}\left(I(q^k) + I(n)\right) < {q^{2k}}I(q^k) + {n^2}I(n).$$
Thus,
$$n\left[q^k - n\right]I(n) < {q^k}\left[q^k - n\right]I(q^k).$$
If $q^k < n$, then $q^k - n < 0$.  Hence,
$$q^k < n \Longrightarrow {q^k}I(q^k) < nI(n) \Longrightarrow \sigma(q^k) < \sigma(n).$$
If $n < q^k$, then $0 < q^k - n$.  Hence,
$$n < q^k \Longrightarrow nI(n) < {q^k}I(q^k) \Longrightarrow \sigma(n) < \sigma(q^k).$$
Consequently, we have
$$q^k < n \Longleftrightarrow \sigma(q^k) < \sigma(n),$$
as desired.
\end{proof}

\begin{lemma}\label{Lemma9}
Let $N = {q^k}{n^2}$ be an odd perfect number given in Eulerian form.  If
$$\frac{\sigma(q^k)}{n} + \frac{\sigma(n)}{q^k} < I(q^k) + I(n),$$
then
$$q^k < n \Longleftrightarrow \sigma(n) < \sigma(q^k).$$
\end{lemma}

\begin{proof}
The proof of Lemma \ref{Lemma9} is very similar to the proof of Lemma \ref{Lemma8}.
\end{proof}

Now, assume that $$\frac{\sigma(q^k)}{n} + \frac{\sigma(n)}{q^k} < I(q^k) + I(n).$$ 
Consider the conclusion of the implication in Lemma \ref{Lemma9} in light of the result $I(q^k) < I(n)$:
$$q^k < n \Longleftrightarrow \sigma(n) < \sigma(q^k).$$
If $q^k < n$, then since $I(q^k) < I(n)$ implies that
$$\frac{\sigma(q^k)}{\sigma(n)} < \frac{q^k}{n},$$
we have
$$\frac{\sigma(q^k)}{\sigma(n)} < \frac{q^k}{n} < 1,$$
which further implies that $\sigma(q^k) < \sigma(n)$.  This contradicts Lemma \ref{Lemma9}.
Similarly, if $\sigma(n) < \sigma(q^k)$, then
$$1 < \frac{\sigma(q^k)}{\sigma(n)} < \frac{q^k}{n},$$
from which it follows that $n < q^k$.  Again, this contradicts Lemma \ref{Lemma9}.
Hence, we know that
$$n < q^k < \sigma(q^k) < \sigma(n)$$
must hold, under the given assumption.  Assuming Brown's proof for $q^k < n$ is completed, this case is ruled out.
Consequently, the inequality
$$\frac{\sigma(q^k)}{n} + \frac{\sigma(n)}{q^k} < I(q^k) + I(n)$$
cannot be true.
Therefore, the reverse inequality
$$I(q^k) + I(n) \leq \frac{\sigma(q^k)}{n} + \frac{\sigma(n)}{q^k}$$
must be true.

It remains to consider the case when
$$I(q^k) + I(n) = \frac{\sigma(q^k)}{n} + \frac{\sigma(n)}{q^k}.$$
Notice that this is true if and only if
$$\sigma(q^k) = \sigma(n),$$
(because $q^k \neq n$). Thus, since $I(q^k) < I(n)$, this implies that $n < q^k$.  Again, assuming Brown's proof for $q^k < n$ is completed, this case is ruled out.

In other words (by Lemma \ref{Lemma8}), we have Theorem \ref{Theorem1} (and the corollary that follows).

\begin{corollary}\label{Corollary1}
If $N = {q^k}{n^2}$ is an odd perfect number given in Eulerian form, then the following biconditional is true:
$$q^k < n \Longleftrightarrow \frac{\sigma(q^k)}{n} < \frac{\sigma(n)}{q^k}.$$
\end{corollary}

We now give another condition that is equivalent to the author's conjecture (mentioned in the introduction).

\begin{theorem}\label{Theorem2}
If $N = {q^k}{n^2}$ is an odd perfect number given in Eulerian form, then the following biconditional is true:
$$\frac{\sigma(q^k)}{n} < \frac{\sigma(n)}{q^k} \Longleftrightarrow \frac{q^k}{n} + \frac{n}{q^k} < \frac{\sigma(q^k)}{\sigma(n)} + \frac{\sigma(n)}{\sigma(q^k)}.$$
\end{theorem}

\begin{proof}
Let $N$ be an odd perfect number given in Eulerian form.  Then $N = {q^k}{n^2}$ where $q \equiv k \equiv 1 \pmod 4$ and $\gcd(q, n) = 1$.

First, we show that
$$\frac{\sigma(q^k)}{n} < \frac{\sigma(n)}{q^k}$$
implies
$$\frac{q^k}{n} + \frac{n}{q^k} < \frac{\sigma(q^k)}{\sigma(n)} + \frac{\sigma(n)}{\sigma(q^k)}.$$
Since $I(q^k) < I(n)$, we have that
$$\frac{\sigma(q^k)}{\sigma(n)} < \frac{q^k}{n}.$$
On the other hand, the inequality
$$\frac{\sigma(q^k)}{n} < \frac{\sigma(n)}{q^k}$$
gives us that
$$\frac{\sigma(q^k)}{\sigma(n)} < \frac{n}{q^k}.$$
This in turn implies that
$$\frac{q^k}{n} < \frac{\sigma(n)}{\sigma(q^k)}.$$
Putting these inequalities together, we have the series
$$\frac{\sigma(q^k)}{\sigma(n)} < \frac{q^k}{n} < \frac{\sigma(n)}{\sigma(q^k)}.$$
Now consider the product
$$\left(\frac{\sigma(q^k)}{\sigma(n)} - \frac{q^k}{n}\right)\left(\frac{\sigma(n)}{\sigma(q^k)} - \frac{q^k}{n}\right).$$
This product is negative.  Consequently we have
$$\left(\frac{\sigma(q^k)}{\sigma(n)}\right)\left(\frac{\sigma(n)}{\sigma(q^k)}\right) - \left(\frac{q^k}{n}\right)\left(\frac{\sigma(q^k)}{\sigma(n)} + \frac{\sigma(n)}{\sigma(q^k)}\right) + {\left(\frac{q^k}{n}\right)}^2 < 0,$$
from which it follows that
$$1 + {\left(\frac{q^k}{n}\right)}^2 < \left(\frac{q^k}{n}\right)\left(\frac{\sigma(q^k)}{\sigma(n)} + \frac{\sigma(n)}{\sigma(q^k)}\right).$$
Therefore, we obtain
$$\frac{n}{q^k} + \frac{q^k}{n} < \frac{\sigma(q^k)}{\sigma(n)} + \frac{\sigma(n)}{\sigma(q^k)}$$
as desired.

Next, assume that $$\frac{\sigma(n)}{q^k} < \frac{\sigma(q^k)}{n}.$$
Since $I(q^k) < I(n)$, we obtain $$\frac{n}{q^k} < \frac{\sigma(q^k)}{\sigma(n)} < \frac{q^k}{n}.$$
Now consider the product
$$\left(\frac{n}{q^k} - \frac{\sigma(q^k)}{\sigma(n)}\right)\left(\frac{q^k}{n} - \frac{\sigma(q^k)}{\sigma(n)}\right).$$
This product is negative.  Therefore, we obtain
$$\left(\frac{n}{q^k}\right)\left(\frac{q^k}{n}\right) - \left(\frac{\sigma(q^k)}{\sigma(n)}\right)\left(\frac{n}{q^k} + \frac{q^k}{n}\right) + {\left(\frac{\sigma(q^k)}{\sigma(n)}\right)}^2 < 0,$$
from which we get
$$1 + {\left(\frac{\sigma(q^k)}{\sigma(n)}\right)}^2 < \left(\frac{\sigma(q^k)}{\sigma(n)}\right)\left(\frac{n}{q^k} + \frac{q^k}{n}\right).$$
Consequently, we have
$$\frac{\sigma(n)}{\sigma(q^k)} + \frac{\sigma(q^k)}{\sigma(n)} < \frac{n}{q^k} + \frac{q^k}{n}.$$
Together with the result in the previous paragraph, this shows that
$$\frac{\sigma(q^k)}{n} < \frac{\sigma(n)}{q^k}$$
is equivalent to 
$$\frac{q^k}{n} + \frac{n}{q^k} < \frac{\sigma(q^k)}{\sigma(n)} + \frac{\sigma(n)}{\sigma(q^k)}.$$
\end{proof}

\begin{remark}\label{Remark8}
Let $N = {q^k}{n^2}$ be an odd perfect number given in Eulerian form. \\
 
Note that, in general, it is true that
$$\frac{\sigma(q^k)}{\sigma(n)} + \frac{\sigma(n)}{\sigma(q^k)} < \frac{\sigma(q^k)}{n} + \frac{\sigma(n)}{q^k},$$
and
$$\frac{q^k}{n} + \frac{n}{q^k} < \frac{\sigma(q^k)}{n} + \frac{\sigma(n)}{q^k}.$$
Therefore,
$$\frac{\sigma(q^k)}{n} < \frac{\sigma(n)}{q^k}$$
is equivalent to
$$\frac{q^k}{n} + \frac{n}{q^k} < \frac{\sigma(q^k)}{\sigma(n)} + \frac{\sigma(n)}{\sigma(q^k)} < \frac{\sigma(q^k)}{n} + \frac{\sigma(n)}{q^k},$$
while
$$\frac{\sigma(n)}{q^k} < \frac{\sigma(q^k)}{n}$$
is equivalent to
$$\frac{\sigma(q^k)}{\sigma(n)} + \frac{\sigma(n)}{\sigma(q^k)} < \frac{q^k}{n} + \frac{n}{q^k} < \frac{\sigma(q^k)}{n} + \frac{\sigma(n)}{q^k}.$$
\end{remark}

At this point, we dispose of the following lemma:

\begin{lemma}\label{Lemma10}
Let $N = {q^k}{n^2}$ be an odd perfect number given in Eulerian form.  Then at least one of the following sets of inequalities is true:
\begin{itemize}
{
\item{$\bf{A}: q^k < \sigma(q^k) < n < \sigma(n)$}
\item{$\bf{B}: q^k < n < \sigma(q^k) < \sigma(n)$}
\item{$\bf{C}: n < q^k < \sigma(n) < \sigma(q^k)$}
\item{$\bf{D}: n < \sigma(n) \leq q^k < \sigma(q^k)$}
}
\end{itemize}
\end{lemma}

Lemma \ref{Lemma10} is proved by listing all possible permutations of the set
$$\left\{q^k, n, \sigma(q^k), \sigma(n)\right\}$$
and then using Theorem \ref{Theorem1}.

Note that Brown's result that $q^k < n$, when completed, would rule out cases $\bf{C}$ and $\bf{D}$ in Lemma \ref{Lemma10}.  Also, notice that by assuming $k=1$, case $\bf{B}$ is also ruled out.

Consequently, we have the following theorem.

\begin{theorem}\label{Theorem3}
Let $N = {q^k}{n^2}$ be an odd perfect number given in Eulerian form.  If $k=1$, then $\sigma(q^k)<n$.
\end{theorem}  

As a corollary, by the contrapositive to Theorem \ref{Theorem3}, we have:

\begin{corollary}\label{Corollary2}
Let $N = {q^k}{n^2}$ be an odd perfect number given in Eulerian form.  If $n<\sigma(q^k)$, then $k>1$.
\end{corollary}

\begin{remark}\label{Remark9}
If one could show the biconditional
$$n<q^{k+1} \iff n<\sigma(q^k),$$
then one would be able to show that
$$n<q^{k+1} \implies k > 1.$$
By the contrapositive, one would then have
$$k=1 \implies q^{k+1} < n \implies q^2 < n.$$
However, we know that
$$n<q^2 \implies k=1.$$
Consequently,
$$n<q^2 \implies k=1 \implies q^2 < n$$
which proves that $q^2 < n$, strengthening Brown's result.
\end{remark}

\section{Final Analysis of the New Results}\label{Section4}
The new results presented in this article seem to imply the following conjecture (see \cite{Dris2}).

\begin{conjecture}\label{Conjecture1}
Let $N = {q^k}{n^2}$ be an odd perfect number given in Eulerian form.  Then the Descartes-Frenicle-Sorli conjecture is false.  (That is, $k > 1$ must hold.)
\end{conjecture}

\begin{remark}\label{Remark10}
Notice how all of the implications in the Lemmas \ref{Lemma1}, \ref{Lemma2} and \ref{Lemma3} in Section \ref{Section2} become vacuously true, given Brown's result that \\$q < n$.  Also, notice that, in Section \ref{Section3}, we could specialize Theorem \ref{Theorem1} (and its consequences) to the case $k=1$ and still get the same results, as follows:
$$q < n \iff \sigma(q) < \sigma(n) \iff \frac{\sigma(q)}{n} < \frac{\sigma(n)}{q}.$$
\end{remark}

\section{Conclusion}\label{Section5}
An improvement to the currently known upper bound of $I(n) < 2$ will be considered a major breakthrough.  In the sequel (\url{http://arxiv.org/abs/1303.2329}), a viable approach towards improving the inequality $I(n) < 2$ will be presented, which may necessitate the use of ideas from the paper \cite{Ward}.

\section{Acknowledgments}
The author sincerely thanks the anonymous referee(s) who have made several corrections and suggestions, which helped in improving the style of the paper.  The author would like to thank Pascal Ochem for communicating the sharper lower bound for $I(n)$.  The author also wishes to thank Carl Pomerance for pointing out the relevance of the paper \cite{AcquaahKonyagin}.  The author also expresses his gratitude to Peter Acquaah for helpful e-mail exchanges on the topic.  Lastly, the author expresses his gratitude to an anonymous reader ``Pascal" who pointed out some ``errors'' in an earlier version of \cite{Dris} (see \cite{Dris2}), thus encouraging him to come up with this sequel.

\end{document}